\documentstyle[11pt,seceq,amssymb]{article}

\newcommand{\be}{\begin{equation}}
\newcommand{\ee}{\end{equation}}
\newcommand{\lb}{\label}
\newcommand{\en}{\varepsilon}
\newcommand{\om}{\omega}
\newcommand{\bl}{{\bf \ell}}
\newcommand{\bu}{{\bf u}}
\newcommand{\bx}{{\bf x}}
\newcommand{\bK}{{\bf K}}
\newcommand{\bD}{{\bf D}}
\newcommand{\tor}{{\Bbb T}^2}
\newcommand{\bphi}{{\mbox{\boldmath $\phi$}}}
\newcommand{\bsigma}{{\mbox{\boldmath $\sigma$}}}
\newcommand{\grad}{{\mbox{\boldmath $\nabla$}}}
\newcommand{\bdot}{{\mbox{\boldmath $\cdot$}}}
\newcommand{\btimes}{{\mbox{\boldmath $\times$}}}
\newcommand{\bzed}{{\mbox{\boldmath $0$}}}
\newtheorem{Th}{Theorem}

\newtheorem{Lm}{Lemma}

\newtheorem{Cor}{Corollary}
\newtheorem{Con}{Conjecture}
\newtheorem{Prop}{Proposition}

\begin{document}
\title{Dissipation in Turbulent Solutions of 2-D Euler}
\author{Gregory L. Eyink\\{\em Department of Mathematics}\\
{\em University of Arizona}\\{\em Tucson, AZ 85721}}
\date{ }
\maketitle
\begin{abstract}
We establish local balance equations for smooth functions of the vorticity
in the DiPerna-Majda
weak solutions of 2D incompressible Euler, analogous to the balance proved by
Duchon and Robert
for kinetic energy in 3D. The anomalous term or defect distribution therein
corresponds to the
``enstrophy cascade'' of 2D turbulence. It is used to define a rather
natural notion of
``dissipative Euler solution'' in 2D. However, we show that the DiPerna-Majda
solutions with vorticity in $L^p$ for $p>2$ are conservative and have zero
defect. Instead, we must seek an alternative approach to dissipative
solutions in 2D. If we assume an upper bound on the energy spectrum of 2D
incompressible Navier-Stokes solutions by the Kraichnan-Batchelor $k^{-3}$
spectrum, uniformly for high Reynolds number, then we show that the
zero viscosity limits of the Navier-Stokes solutions exist, with vorticities
in the zero-index Besov space $B^{0,\infty}_2$, and that these give a
weak solution of the 2D incompressible Euler equations.
We conjecture that for this class of weak solutions enstrophy dissipation
may indeed occur, in a sense which is made precise.
\end{abstract}

\newpage

\section{Introduction}

In 2-dimensional turbulence it is the {\it enstrophy}
$\Omega(t):={{1}\over{2}}\|\omega(t)\|_2^2$ that is
expected to cascade to small length-scales, not the energy as in three space
dimensions \cite{Kr67}-\cite{Ba69}.
In a view that goes back to Onsager \cite{On49}, such turbulent cascades are
conjectured to be described, in
the limit of infinite Reynolds number, by singular (or weak) solutions of the
incompressible Euler equations.
More recently, Duchon and Robert \cite{DR00} have shown how Onsager's idea of a
dissipative Euler solution
may be formalized in the three-dimensional case via a local energy balance
relation. It is our purpose here
to similarly formalize the notion of a 2-dimensional dissipative Euler
solution, corresponding to the enstrophy cascade.

We consider weak solutions of the 2D Euler equations in the
vorticity-velocity formulation:
\be \partial_t\om + (\bu\bdot\grad)\om = 0, \lb{euler} \ee
with $\bu = \bK*\omega$ given by the Biot-Savart kernel $\bK$.
We show first that when the vorticity fields $\om(\bx,t)$ are
suitable measurable functions and (\ref{euler}) is interpreted in the
sense of distributions, then a local balance is satisfied
\be \partial_th(\om) + \grad\bdot[\bu h(\om)] = -Z_h(\omega), \lb{h-bal} \ee
for nonnegative, convex functions $h(\omega)$. Of course, (\ref{euler})
formally just expresses the conservation of vorticity along fluid particle
trajectories, so that (\ref{h-bal}) would naively be expected to hold
with $Z_h(\om)\equiv 0.$  Dissipative weak solutions might
be taken to be those for which this distribution is nonnegative:
$Z_h(\omega)\geq 0$. The balance equation (\ref{h-bal}) makes
more precise Polyakov's analogy of the enstrophy cascade with
conservation law anomalies in quantum field-theory (such as the axial
anomaly in QED) \cite{Poly}. The distribution appearing as a sink term
on the right side of (\ref{h-bal}) corresponds closely to such an anomaly.
However, we show under rather general conditions, even weaker than those
in our earlier work \cite{Ey96}, that $Z_h(\omega)\equiv 0$. For example,
we show that the anomaly vanishes for functions $h$ of ``power-$p$ growth''
at large arguments, whenever the initial vorticity satisfies an $L^p$ bound
in space. In particular, this means that {\it enstrophy is conserved
by a 2D Euler solution whenever the enstrophy itself is finite.}
\footnote{This statement has, among other consequences, the implication
that no power-law 2D energy spectrum $E(k)\sim C k^{-n}$ with $n>3$
may be compatible with existence of an enstrophy cascade. The conformal
``solutions'' studied by Polyakov \cite{Poly} that have spectral exponent
$n>3$ therefore cannot exhibit an anomaly in the enstrophy conservation
law, as he has proposed.} This state of affairs presents a striking
contrast with the situation in 3D where, as discussed by Duchon and Robert
\cite{DR00}, energy dissipation is expected to be possible for
incompressible Euler solutions with finite energy.

The above results necessitate an approach to the notion of dissipative
Euler solution in the 2D case which is therefore rather different from
that of Duchon-Robert for 3D. Nevertheless, 2D turbulence theory is still
a useful guide to the correct formulation. Indeed, the above results
are in perfect agreement with the expectations of the classical theories
\cite{Kr67}-\cite{Ba69}, which predict that the small-scale energy
spectrum in the enstrophy cascade range of 2D turbulence shall be
of the form $E(k)\sim C k^{-3}$ (with at most a logarithmic correction).
Hence, the classical theories of 2D turbulence predict an infinite total
enstrophy but a finite spectral flux of enstrophy. In the following
we shall formalize this notion for an appropriate class of weak Euler
solutions in 2D. We define as dissipative those solutions which have
a nonnegative flux of enstrophy (possibly zero or infinite)
asymptotically to infinitely high wavenumbers. The relevant
solutions must, however, have vorticity fields which exist only as
distributions and not as ordinary (measurable) functions. We show that
such solutions of 2D Euler equations with a Kraichnan-Batchelor $k^{-3}$
energy spectrum are obtained as zero-viscosity limits of the Leray
solutions of 2D Navier-Stokes, whenever upper bounds on the spectrum
by the Kraichnan-Batchelor prediction hold uniformly in the viscosity.
We then show that {\it the notion of enstrophy flux is well-defined for
such distributional solutions of 2D Euler, although the enstrophy
itself may diverge.} Our natural conjecture is that the flux is
asymptotically nonnegative at small length-scales for all such
``viscosity solutions'' and, for suitable initial data, even
strictly positive.

Our main results are stated as three Theorems in the following
Section 2, where their content is further discussed in detail.
The proofs of the Theorems are outlined in the final Section 3.

\newpage

\section{Statement of Results}

Before stating precisely our theorems, it will help to motivate the statements
(and the proofs) to give a brief,
heuristic argument for the existence of the enstrophy cascade. In \cite{Ey96}
we considered a ``filtered'' form of the 2-D Euler equations
(see also \cite{C-W}):
\be \partial_t \om_\en + \grad\bdot[\bu_\en \om_\en + \bsigma_\en] = 0,
\lb{filt-E} \ee
where $\omega_\en = \varphi_\en * \omega$ for a smooth mollifier $\varphi$,
$\varphi_\en(\bx)= \en^{-2}\varphi(\en^{-1}\bx),$
and $\bsigma_\en = (\bu \om)_\en -\bu_\en\om_\en$. The new term $\bsigma_\en$
represents a turbulent spatial transport of vorticity
due to the eliminated small-scales. It is straightforward to show that the
balance holds that
\be \partial_th(\om_\en) + \grad\bdot[\bu_\en h(\om_\en) +
h'(\om_\en)\bsigma_\en] =
h''(\om_\en)\grad\om_\en\bdot \bsigma_\en. \lb{filt-bal} \ee
The term $Z_{h,\en}(\omega):= -h''(\om_\en)\grad\om_\en\bdot \bsigma_\en$
represents a transfer of $h$-stuff from length-scales
$>\en$ to smaller scales. Based upon the notion of ``UV-locality of
interactions'', a natural approximation is to take $\bsigma_\en \approx
({\rm const.})[(\bu_\en\om_\en)_\en-\bu_\en\om_\en]$ and then to Taylor expand
to leading non-vanishing order to obtain
\be \bsigma_\en \approx C \en^2 \bD_\en \bdot \grad\om_\en. \lb{nlin-mod} \ee
Here $\bD_\en$ is the filtered velocity-gradient tensor $D_{ij}=\partial
u_i/\partial x_j$; also, a spherically symmetric mollifier
has been assumed. The first of our approximations is analogous to the
``similarity model'' employed by engineers in large-eddy simulation
of three-dimensional turbulence and the second to its further simplification,
the ``nonlinear model'' \cite{MK00}. The matrix
$\bD_\en$ is traceless and has, in vortical regions of the flow, a pair of
imaginary eigenvalues and, in strain-dominated
regions, two real eigenvalues of equal magnitude $S_\en$ but opposite signs. It
stands to reason that, in the latter straining
regions, the compression of vorticity level sets will tend to align the
direction of the vorticity gradient $\grad\om_\en$
with the eigendirection of $\bD_\en$ corresponding to the negative eigenvalue.
Indeed, such alignment has been observed
in simulations to hold (for the unfiltered quantities) with a high probability
\cite{BMPS}. Assuming it to hold exactly, we find that
\be \bsigma_\en \approx -C\en^2 S_\en \grad\om_\en.  \lb{eddy-visc} \ee
This is precisely an eddy-viscosity model, with effective viscosity
$\nu_\en=C\en^2 S_\en$ at scale $\en$. It leads to
an effective dissipation $Z_{\en}(\om)\approx -\nu_\en |\grad\om_\en|^2.$ If
the vorticity field is H\"{o}lder continuous
with exponent $s$, $\omega\in C^s$, then $\grad\om_\en\sim \en^{s-1}$ for small
$\en$ and $S_\en\sim S$ independent of $\en$.
In that case, $Z_{\en}(\om)\sim \en^{2s}$ for $\en\rightarrow 0$, so that we
expect an asymptotic enstrophy cascade
only when $s=0$. This is precisely the ``mean-field'' scaling exponent in the
Batchelor-Kraichnan theory \cite{Kr67}-\cite{Ba69}.

We now state our main theorems:

Our first theorem establishes the local vorticity balance equations for the
weak Euler solutions constructed by
DiPerna and Majda for initial data $\omega_0\in L^p, p>1$ \cite{DP-M}. Although
they considered solutions in the whole
plane ${\Bbb R}^2$, we shall restrict attention for simplicity to solutions on
the 2-D torus $\tor$. DiPerna and Majda
also established existence of weak solutions in the velocity-pressure
formulation, but it is not hard to show that, for $p>4/3,$
the associated vorticity field in their solution also satisfies the weak
vorticity-velocity equations (see below).
In fact, the only property of the DiPerna-Majda solution that we will employ in
our proof is that $\omega\in
L^\infty([0,T],L^p(\tor))$ and our theorem would apply to any other such
solutions as well. To state our theorem,
we must introduce an appropriate class of differentiable functions
\be {\cal H}_p := \left\{ h |\,\,\,h\in C^1({\Bbb R}),
\,\,\,|h'(\omega)|\leq C(1+|\omega|^{p-1})\,\,{\rm for}\,\,{\rm some}
                      \,\, C>0 \right\} \lb{Cp1} \ee
which have at most $L^p$-growth. We then have the following:
\begin{Th}
If $\om\in L^\infty([0,T],L^p(\tor))$ and the associated $\bu=\bK*\om$ for
$p>4/3$ are a weak solution
of 2-D incompressible Euler in the vorticity-velocity formulation, then for
$h\in {\cal H}_r\bigcap C^2$, with $r={{3}\over{2}}p-1$ for
${{4}\over{3}}<p < 2, r<p$ for $p=2,$ and $r=p$ for $p>2,$ the balance
(\ref{h-bal}) holds
$$ \partial_t h(\om) + \grad\bdot[\bu h(\om)] = - Z_h(\om) $$
in the sense of distributions. The righthand side is given by the
distributional limit
\be  Z_h(\om) = \lim_{\en\rightarrow 0}
-h''(\om_\en)\grad\om_\en\bdot\bsigma_\en \lb{Zh} \ee
which exists for any choice of mollifier $\varphi$ which is $C^\infty$,
nonnegative, and compactly supported, with unit integral,
and it is independent of that choice. For the special case of the enstrophy
integral, $h(\omega)={{1}\over{2}}|\om|^2,$ when $p>2$,
we write simply $Z(\om)=  Z_h(\om)$. In that case, there is the alternative
expression:
\be Z(\omega) = \lim_{\en\rightarrow 0} {{1}\over{4}}\int d^2\bl
\,\,\grad\varphi_\en(\bl)\bdot\Delta_\bl\bu|\Delta_\bl\om|^2 \lb{4/5law} \ee
where $\Delta_\bl\om(\bx,t)=\om(\bx+\bl,t)-\om(\bx,t),$ likewise for
$\Delta_\bl\bu$, and $\varphi$ is further restricted to be an
even function of its argument.
\end{Th}
Note that, formally, $ Z_h(\om) = h''(\om) Z(\om)$, so the fluxes of general
convex functions are, in some sense, proportional
to the enstrophy flux with a nonnegative factor. The last expression
(\ref{4/5law}) for the enstrophy flux has a nice
interpretation as a local, non-ensemble-averaged form of the ``-2 law'' for the
direct cascade, in its form applicable
without isotropy (see \cite{Ey96}, Appendix B). Thus, the defect distribution
in the vorticity balance equations has an exact connection with the
enstrophy cascade in 2D turbulence theory.

However, we next show that this distribution is, in fact, zero for the
DiPerna-Majda weak solutions, which therefore conserve the integral
\be                I_h(t) = \int d^2\bx \,\, h(\om(\bx,t))     \lb{Ih} \ee
for all $h$ of suitable growth:
\begin{Th}
If $\om\in L^\infty(0,T;L^p(\tor))$ is a DiPerna-Majda weak Euler solution
for $p\geq 2$, then
\be \partial_t h(\om) + \grad\bdot[\bu h(\om)] =0 \lb{loc-cons} \ee
in distribution sense for all $h\in {\cal H}_r$ with $r=p$ when $p>2$ and
for any $r<p$ when $p=2$.
\end{Th}
In \cite{Ey96} it was proved that such a conservation statement holds
for $\om\in L^p(0,T;B^{s,\infty}_p(\tor))$ for $s>0,p\geq 3$
where $B_p^{s,\infty}({\Bbb T}^2)$ is the standard Besov space
of functions in $L^p({\Bbb T}^2)$ which are H\"{o}lder of index
$s$ in the $L^p$-mean sense \cite{Treib}. That theorem was thus analogous to
the Besov-space improvement of Onsager's original conservation
result for 3D, which was proved by Constantin, E, and Titi \cite{CET}.
We now see that the smoothness assumed in \cite{Ey96} was unnecessary and
that simple $L^p$ bounds alone are sufficient for conservation. Essentially
the same result was already obtained by P.-L. Lions in \cite{Lions96},
Section 4.1, based upon his earlier work with R. J. DiPerna \cite{DP-L}.
He showed there that the DiPerna-Majda solutions with $p>2$ are
``renormalized solutions'' in the sense of DiPerna-Lions \cite{DP-L},
which amounts to the requirement that (\ref{loc-cons}) hold. In fact,
global conservation
\be \int_{\tor} d^2\bx\,\,h(\om(\bx,t)) =
 \int_{\tor} d^2\bx\,\,h(\om_0(\bx)), \,\,\,\,\,t>0 \lb{glob-cons} \ee
is shown in \cite{Lions96} to hold  for all $h\in {\cal H}_p$ even when
$p=2$, just as in the proof of Theorem II.2 and equation (26) in
DiPerna-Lions \cite{DP-L}. \footnote{The same remark was made in a recent
preprint of E and Vanden-Eijnden \cite{E-vdE2}.} In particular, taking
$h(\om)={{1}\over{2}}|\om|^2$, a remarkable statement is true that
{\it enstrophy dissipation is not possible for any 2D Euler solutions
with finite enstrophy}. We conclude more generally that the DiPerna-Majda
weak solutions are not relevant to the problem of constructing dissipative
Euler solutions. In the language of turbulence theory, they do not support
enstrophy cascades over infinitely-long ranges of wavenumber.

The conservation properties of the DiPerna-Majda solutions for $p>2$
have an intuitive explanation. It has been noted recently that breakdown of
uniqueness of Lagrangian particle trajectories in H\"{o}lder but non-Lipschitz
flows can be a mechanism for the anomalous dissipation
of the analogous integrals as (\ref{Ih}) for passive scalars
\cite{GKB}-\cite{E-vdE}.
For the 3D problem, Shnirelman has found a weak solution which dissipates
energy globally, by constructing a generalized
flow with random Lagrangian trajectories \cite{Shnir}. In the case of
the Yudovich solutions of 2D Euler with $\om\in L^\infty(\tor)$
\cite{Yud}, it has long been known that they are conservative
precisely because the corresponding velocity field is log-Lipschitz and the
Lagrangian flow maps $X_t$ are unique, volume-preserving homeomorphisms.
Therefore, the Yudovich solution is given simply by
$\om(\bx,t)=\om_0(X_{-t}(\bx))$ in terms of the inverse-Lagrangian map.
All of the integrals $I_h(t)$ in (\ref{Ih}) are then trivially time-invariant.
DiPerna and Lions in their paper \cite{DP-L} show that there are
likewise unique Lagrangian flow maps $X_t(x)$ with $X\in C(0,T;L^p(\tor))$
whenever $\bu\in L^1(0,T;W^{1,p}(\tor))$ for $p\geq 1$ and that these maps
preserve Lebesgue measure when $\grad\bdot\bu=0$. The ``renormalized''
solutions of the linear advection equation constructed by DiPerna-Lions
are shown to have precisely the form $\om(\bx,t)=\om_0(X_{-t}(\bx))$.
While it is not true in general that the distributional solutions of 2D Euler
in the sense considered here are renormalized solutions, Theorem 2, as we have
observed above, shows that this is so for the DiPerna-Majda solutions when
$p>2$. Hence, the conservation properties of these solutions are again
connected with the uniqueness of Lagrangian particle trajectories.

The above results are negative---in the sense that they imply a lack of
enstrophy
dissipation---but we wish to emphasize that they are fully consistent with the
expectations of 2D turbulence theory. In fact, the Navier-Stokes solutions
exhibiting an enstrophy cascade are expected to have the Batchelor-Kraichnan
energy spectrum
\be E(k,t) \sim C\eta^{2/3}(t)k^{-3}, \lb{KrBa-spec} \ee
where $\eta(t)$ is the enstrophy dissipation rate per volume \cite{Kr67,Ba69}.
This spectrum should hold at high wavenumbers $k\gg k_0(t)$, the wavenumber
of peak enstrophy, up to a wavenumber $k_d(t)= \nu^{-1/2}\eta^{1/6}(t)$,
at which the dissipation by viscosity $\nu$ becomes relevant. Equivalently,
the enstrophy spectrum predicted by Batchelor-Kraichnan theory is
\be \Omega(k,t) \sim C\eta^{2/3}(t)k^{-1}, \lb{KrBa-enst-spec} \ee
for $k_0(t)\ll k\ll k_d(t)$. In the limit as $\nu\rightarrow 0$
this spectrum extends all the way to $+\infty$ and its integral diverges,
implying an infinite total enstrophy. As we have seen, this is rigorously
required to have limiting Euler solutions which can dissipate enstrophy.

In fact, velocity fields $\bu(t)$ with the Batchelor-Kraichnan spectrum
(\ref{KrBa-spec}) for all $k\gg k_0(t)$, when that spectrum is interpreted
in a suitable sense, must consist of $\bu(t)\in B^{1,\infty}_2(\tor),$ with
corresponding vorticity $\omega(t)\in B^{0,\infty}_2(\tor)$, the Besov space of
zero index. The definition of spectrum which is relevant is a
``Littlewood-Paley
spectrum''which was earlier used by P. Constantin in \cite{Const-spect}
to prove a rigorous upper bound. This spectrum is defined in terms of the
Littlewood-Paley decomposition of the velocity $\bu(t)=\sum_{N=0}^\infty
\bu_N(t)$ with $\bu_N(t)=\psi_N*\bu(t)$, for a smooth partition of unity
in wavenumber space
\be \widehat{\psi}_0(k)+\sum_{N=1}^\infty \widehat{\psi}_N(k)=1
                    \lb{LP-part} \ee
where ${\rm supp}(\widehat{\psi}_N)\subset [2^{N-1},2^{N+1}]$ for $N\geq 1$
and ${\rm supp}(\widehat{\psi}_0)\subset [0,2]$. The Littlewood-Paley
spectrum is then defined by
\be E_{LP}(k,t):= k^{-1}\|\bu_N(t)\|_{L^2}^2 \lb{LP-spec} \ee
for $k\in [2^N,2^{N+1})$. With this definition it is not hard to see
that $\bu(t)\in B^{1,\infty}_2(\tor)$ precisely when the spectrum satisfies
a bound of the form $E_{LP}(k,t)=O(k^{-3})$. In fact, the Littlewood-Paley
criterion for $f\in B^{s,\infty}_p(\tor)$ with $-\infty<s<\infty$ and $p>0$
is just that
\be   \|f\|_{B^{s,\infty}_p} :=  \sup_{N\geq 0} 2^{sN}\|f_N\|_{L^p}
       \lb{Besov-norm} \ee
be finite, and for $p\geq 1$ this is a norm making $B^{s,\infty}_p(\tor)$
into a Banach space. See \cite{Treib}, Sections 2.3.1-3. Note we assume
only big-$O$ bounds on the spectrum and not a power-law scaling. In fact,
not long after his first paper on 2D turbulence, Kraichnan argued that
there should be a logarithmic correction, $E(k,t)\sim C\eta^{2/3}(t)
k^{-3}[\ln(k/k_0)]^{-{{1}\over{3}}},$ where $k_0$ is the lower end of
the enstrophy cascade range \cite{Kr71}. In any case, it is still generally
believed that the true energy spectrum must be bounded above by the form
(\ref{KrBa-spec}) at high Reynolds number. In that case, we see that
the vorticity field $\omega(t)$ has the Besov index $s=0$ but not
necessarily any larger index. As we have already remarked, $\omega(t)\in
B^{s,\infty}_p(\tor)$ with $p>2$ and $s>0$ could not be consistent with
a non-vanishing enstrophy dissipation.

It is still an open question whether solutions of 2D Euler equations
exist with velocities and vorticities in such Besov spaces and, if
so, whether they dissipate enstrophy in a suitable sense. We shall
prove here a few simple results in this direction and, in particular,
advance our main conjecture.

We show first that an upper bound on the energy spectrum of the solutions
of the 2D Navier-Stokes solutions $\bu^\nu(t)$ by the Batchelor-Kraichnan
spectrum, when that bound is uniform in the viscosity, implies the existence
of 2D Euler solutions $\bu(t)$ in the appropriate Besov spaces. In
\cite{Const-spect} Constantin proved that the long-time average energy
spectrum,
$\overline{E}_{LP}^\nu(k)=\limsup_{T\rightarrow\infty}{{1}\over{T}}\int_0^T
dt\,\,E_{LP}^\nu(k,t)$ of the 2D Navier-Stokes solutions satisfies a bound
of the form
\be \overline{E}_{LP}^\nu(k) \leq C\gamma^2 k^{-3}
                           \left({{k_d}\over{k}}\right)^6, \lb{const-bd} \ee
where $\gamma=\|\grad\bu^\nu\|_{L^\infty}$. This upper bound is much larger
than
the Kraichnan-Batchelor spectrum (\ref{KrBa-spec}) over most of the range
$k_0<k<k_d$,
but becomes comparable at the upper limit. We are going to assume here that
something stronger is true of the 2D Navier-Stokes solutions, namely, for
some $T>0$:
\be  \sup_{\nu>0}{{1}\over{T}}\int_0^T dt\,\,\sup_{k>k_0} k^3E_{LP}^\nu(k,t)
     <\infty                    \lb{u-Hyp} \ee
Note that $C^\nu(t)=\sup_{k>k_0}k^3E_{LP}^\nu(k,t)/\eta^{2/3}$ is an
instantaneous
(worst) value of the Kraichnan-Batchelor constant and (\ref{u-Hyp}) is a bound
on its time-average. The most important aspect of this estimate---in contrast
to what is so far proved, equation (\ref{const-bd})---is its uniformity for
small
viscosity $\nu>0$. Our main hypothesis (\ref{u-Hyp}) is equivalent to
\be \sup_{\nu>0}\|\bu^\nu\|_{L^2(0,T;B^{1,\infty}_2(\tor))}<\infty
    \lb{u-besov} \ee
Using estimates for singular integral operators, this may also be expressed
equivalently
in terms of the vorticity $\om^\nu =\grad\btimes\bu^\nu$, as
\be \sup_{\nu>0}\|\om^\nu\|_{L^2(0,T;B^{0,\infty}_2(\tor))}<\infty
     \lb{om-besov} \ee
This latter estimate could be stated in terms of the Littlewood-Paley enstrophy
spectrum $\Omega_{LP}^\nu(k,t):=k^{-1}\|\om_N^\nu(t)\|_{L^2}^2$ for
$k\in [2^N,2^{N+1})$, that
\be \sup_{\nu>0} {{1}\over{T}}\int_0^T dt\,\,\sup_{k>k_0}
    k\Omega_{LP}^\nu(k,t) <\infty, \lb{om-Hyp} \ee
which is entirely equivalent to the initial hypothesis (\ref{u-Hyp}).

We now state our third main result:

\begin{Th}
Let $\bu^\nu$ be the solution of the $2D$ incompressible Navier-Stokes equation
for initial data $\bu_0\in B^{1,\infty}_2(\tor)$ and viscosity $\nu>0$. Let
$\om^\nu=\grad\btimes\bu^\nu.$ Assume that (\ref{u-besov}) holds for these
solutions.
Then, there exists a $\bu\in L^2(0,T;B^{1,\infty}_2(\tor))\cap
{\rm  Lip}(0,T;H^{-L}(\tor))$ which is a weak solution of the 2D incompressible
Euler equations in the velocity-pressure formulation, and for which, with
$\om=\grad\btimes\bu$,
\be \om^\nu\rightharpoonup \om \,\,\,\, weak-* \,\,\,\, {\rm in} \,\,\,\,
      L^2(0,T;B^{0,\infty}_2(\tor)), \lb{nu-om-weak-*} \ee
\be \om^\nu\rightarrow \om \,\,\,\, strong \,\,\,\, in \,\,\,\,
      L^2(0,T;W^{-s,q}(\tor)), \lb{nu-om-strong} \ee
for some $q>2$ and $s>1-{{2}\over{q}}$, in the limit as $\nu\rightarrow 0$.
Furthermore, $\om$ and $\bu=\bK*\om$ are a weak solution of the 2D
incompressible
Euler equations in the vorticity-velocity formulation, in the sense that
\be \langle (\partial_t+\bu\bdot\grad)\psi,\om\rangle=0 \lb{Bes-2DE} \ee
for all $\psi\in C^\infty_0([0,T]\times\tor)$ and the expression
$\langle (\partial_t+\bu\bdot\grad)\psi,\om\rangle$ is defined as the
evaluation of a continuous linear functional on the element
$\om$ of the Banach space $L^2(0,T;B^{0,\infty}_2(\tor))$.
\end{Th}

\noindent This theorem essentially just states that the estimate
(\ref{u-besov})
provides enough compactness to take limits along subsequences. Obviously, the
hard problem is to prove that a bound such as (\ref{u-besov}), as expected from
2D turbulence theory, really does hold. The theorem could be stated in a
somewhat more general form, with the results on weak solutions in the
velocity-pressure formulation remaining true for any $p>1$ replacing
$p=2$, if a corresponding replacement is made in the estimate (\ref{u-besov}).
Likewise, the results on weak solutions in the vorticity-velocity formulation
will remain true for any $p>4/3$ replacing $p=2$.

Our interest in this class of solutions is that they seem compatible with a
finite rate of enstrophy dissipation in the inviscid limit. However, the
very notion of ``dissipative solution'' must be reformulated. Local functions
of the vorticity, of the form $h(\om(\bx,t))$, do not need even to exist,
since now the vorticity $\om$ is only a distribution and not necessarily
a measurable function. Thus, a balance equation such as (\ref{h-bal}) that
we proved in Theorem 1 for DiPerna-Majda solutions is not even well-defined
for the class of solutions considered here. However, the balance equations
for the mollified vorticity in (\ref{filt-bal}), namely,
\be \partial_th(\om_\en) + \grad\bdot[\bu_\en h(\om_\en) +
h'(\om_\en)\bsigma_\en] = -Z_{h,\en}(\om) \lb{filt-bal-rep} \ee
with $Z_{h,\en}(\omega)= -h''(\om_\en)\grad\om_\en\bdot \bsigma_\en,$
are still perfectly well-defined. The term $Z_{h,\en}(\omega)$ which appears
as sink on the righthand side of (\ref{filt-bal-rep}) represents a flux
of $h$ to length-scales $<\en$ and it is expected to be asymptotically
non-negative for small $\en$. In fact, more should be true. A corresponding
balance equation holds for the solutions of the 2D Navier-Stokes solutions
$\om^\nu$, in the form
\be \partial_t h(\om^\nu) + \grad\bdot[h(\om^\nu)\bu^\nu -\nu\grad h(\om^\nu)]
                          = -\nu h''(\om^\nu)|\grad\om^\nu|^2, \lb{NS-h-bal}
\ee
for any $h\in C^2$. Then we expect the following
\begin{Con}
Let $\om^\nu$ be a sequence of solutions of the 2D Navier-Stokes equation
obeying (\ref{u-besov}) and let $\om$ be the limiting 2D Euler solution,
as provided by Theorem 3. Then, for this Euler solution
\be Z_h(\om) = \lim_{\en\rightarrow 0} Z_{h,\en}(\omega) \lb{Z-dist} \ee
exists in the sense of distributions for any $h\in {\cal H}_2$.
Furthermore, the same distribution is obtained by the limit of the
viscous dissipation of the Navier-Stokes solutions:
\be Z_h(\om) = \lim_{\nu\rightarrow 0} \nu h''(\om^\nu)|\grad\om^\nu|^2
    \lb{NS-Z-dist} \ee
for any $h\in {\cal H}_2\bigcap C^2$.In particular, for any such convex
$h$, the distribution $ Z_h(\om)$ is a nonnegative measure. Finally, there
should exist a suitable such 2D Euler solution $\om$ for which
\be Z_h(\om)>0 \lb{strict} \ee
with a {\it strict} inequality, for a convex $h\in {\cal H}_2$.
\end{Con}
The first limit statement in the conjecture may be put another way, which
is perhaps more illuminating. Although the integral $I_h(t)$ may itself be
infinite for the Euler solutions in Theorem 3, it still makes sense to talk
about a finite {\it dissipation rate} for it, defined as $D_h(t):=
\liminf_{\en \rightarrow 0} -{{dI_h^\en}\over{dt}}(t)$, where $I_h^\en(t)$
is the value of the integral for $\om_\en$. The conjecture then states
$D_h(t)=\int_{\tor} d^2\bx\,\,Z_h(\om)(\bx,t)>0$. Note that
for DiPerna-Majda solutions the first limit (\ref{Z-dist})
has been demonstrated in Theorem 1 and it is easy to show for these
solutions that the second limit (\ref{NS-Z-dist}) also holds, using the
same kind of argument as in Proposition 4 of Duchon-Robert \cite{DR00}.
Of course, for DiPerna-Majda solutions with $\om\in L^p$ and $p>2$ the
distribution $Z_h(\om)\equiv 0$ and thus the third statement is false.

We believe that it is necessary to understand solutions of the type
considered in Theorem 3 in order to develop a rigorous mathematical
theory of invariant measures for forced steady-states of 2D Navier-Stokes
in the zero-viscosity limit. As proved in \cite{Ey96}, Section 3.3.4,
the mean enstrophy flux $\langle Z_\en\rangle$ is a positive constant $\eta$,
independent of $\en$, for length-scales $\en\ll \ell_f$, the forcing scale,
and $\en\gg \nu^{1/4}E^{1/4}/\eta^{1/4}$, under the single assumption that
the total mean energy $E$ remains finite in the limit as $\nu\rightarrow 0$.
(This requires adding an additional dissipation at low-wavenumbers
to dispose of the ``condensate'' from the inverse energy cascade: see
\cite{Ey96}).
Thus, in the limit as $\nu\rightarrow 0$, we expect that the realizations
of the ensemble shall be solutions of the (forced) 2D Euler equations
with $Z(\om)>0$. If the statistical energy spectrum has the Batchelor-Kraichnan
form in this limit, then individual realizations of the vorticity satisfy
$\om\in B^{0-,\infty}_2(\tor)$ a.s. It is not hard to prove this fact,
using the methods of \cite{Ey95} (the wavelet characterization of Besov
spaces and the Borel-Cantelli argument of Theorem 4).

To see dissipation in the sense of our Conjecture 1 for the problem of free
decay
of 2D turbulence starting from random initial conditions, one must begin with
initial data which is sufficiently rough. It is well-known that if one starts
with $\om_0\in B^{s,\infty}_p(\tor)$ for $s>0$ at time $t=0$, then the exponent
may
(and generally will) deteriorate exponentially in time: for example, $s(t)=
e^{-C\|\om_0\|_\infty t}s$ if $\om_0\in L^\infty(\tor)\bigcap
B^{s,\infty}_p(\tor),$
but the exponent remains positive \cite{BC94}. Thus, there will be no
dissipation
at any finite time. On the other hand, the deterioration is consistent with the
expectation from 2D turbulence theory that there will be an exponentially
growing
range of scales $\en$ with $\langle Z_\en(t)\rangle \approx \eta(t),$
independent
of $\en$ \cite{Kr75,Kr74}. To see dissipation at finite (or zero) time, one
must
begin with initial data no more regular than $\om_0\in B^{0,\infty}_p(\tor)$
a.s.
for $p\geq 2$. Such initial data could be prepared, for example, by taking an
invariant measure for the driven problem and then turning off the force.
The results of DiPerna and Lions \cite{DP-L} do not rule out dissipation
in this instance, because they require $\bu\in L^1(0,T;W^{1,p}(\tor)$ for
some $p\geq 1$, whereas $W^{1,p}(\tor)= B^{1,p}_p(\tor)\subsetneq
B^{1,\infty}_p(\tor)$. If $\bu(t)\in B^{1,\infty}_p(\tor)$ only,
then examples like that in section IV.2 of \cite{DP-L} show that
uniqueness of the Lagrangian trajectories breaks down and dissipation
(in the sense of non-vanishing enstrophy flux) is possible.

It is natural to expect that 2D Euler solutions which are dissipative
in the proposed sense, i.e. $Z_h(\om)\geq 0$ for convex $h$, must be unique.
Our Conjecture 1 states that ``viscosity solutions'' of 2D Euler equations
are dissipative, so that these must then also be unique. Duchon and Robert
\cite{DR00} have advanced the same idea for the 3D case. There is perhaps
even more reason to believe so in 2D, because there is then an infinity
of convex ``entropies'' $h$. For the problem of scalar conservation laws, such
entropies play a crucial role in establishing uniqueness (e.g. see \cite{Lax}).
However, unlike the scalar case, it is not necessarily true even for smooth
classical solutions of 2D Euler that the dynamics is $L^1$-contractive. In
fact,
for two such solutions $\om_1,\om_2$, ${{d}\over{dt}}\|\om_1(t)-\om_2(t)\|_1
= -2\int_{\om_1=\om_2} {\bf n}_{12}\bdot (\bu_1-\bu_2)\om\,\,ds $ where ${\bf
n}_{12}$ is the
unit vector normal to the curve ${\om_1=\om_2}$ from the region $\om_1>\om_2$
to $\om_2>\om_1$,
and $\om=\om_1=\om_2$. It is precisely the nonlocal relation between $\bu$ and
$\om$ which allows
$\bu_1\neq\bu_2$ where $\om_1=\om_2$. So far, uniqueness
of weak Euler solutions in 2D is established only for the solutions with
$\om\in L^\infty([0,T]\times {\Bbb T}^2)$ constructed by Yudovich \cite{Yud}
and for DiPerna-Majda solutions in $L^p, p\in (1,\infty)$ if also $\om\in BMO$
\cite{Vishik}. It is not known in general whether the DiPerna-Majda solutions
are unique, although for $p>2$ they are ``dissipative Euler solutions'',
in the sense that
\be \partial_t h(\om) + \grad\bdot[\bu h(\om)] \leq  0 \lb{diss-sol} \ee
for all convex $h\in {\cal H}_{p}$. In fact, as noted above, DiPerna-Majda
solutions for $p>2$ are ``renormalized solutions'' as considered by
DiPerna-Lions and satisfy (\ref{diss-sol}) in the degenerate sense with
equality. Yet their uniqueness is an open question.

\newpage

\section{Proofs}

\noindent {\bf 2.1. Proof of Theorem 1}

\noindent We comment first on the validity of the weak vorticity-velocity
equation for the DiPerna-Majda solutions. The condition
$p>4/3$ arises from the requirement that the nonlinear advection term
$\bu\om\in L^1(\tor)$. Since $\bu\in W^{1,p}(\tor)\subset
L^{p'}(\tor)$ for ${{1}\over{p'}}={{1}\over{p}}-{{1}\over{2}}$ by Sobolev
imbedding, one finds that $p'>q$, with $q$  defined by ${{1}\over{q}}
=1-{{1}\over{p}}$, when $p>4/3$. Then $\bu\om\in L^1(\tor)$ follows by
H\"{o}lder inequality. The weak velocity-pressure form
of the Euler equation is that
\be \int d^2\bx \int dt\,\, [\partial_t\bphi\bdot\bu +
\grad\otimes\bphi:\bu\otimes\bu]  = 0 \lb{vel-press} \ee
for any smooth, divergence-free test function $\bphi(\bx,t)$. In particular,
$\bphi = \grad^\perp \psi$ satisfies
these conditions for any smooth $\psi$, where $\grad^\perp$ is the
skew-gradient,
$\partial_i^\perp = \varepsilon_{ij}\partial_j$
with $\varepsilon_{ij}$ the Levi-Civita tensor in 2D. (In fact, by Hodge
theory, any divergence-free vector field $\bphi$ in 2D
can be written in this way.) Substituting $\bphi = \grad^\perp \psi$ into
(\ref{vel-press}) it is easy, using the $L^1$ property of
$\bu\om$ and $\omega= -\grad^\perp\bdot\bu$, to derive the vorticity-velocity
equation by an approximation argument.

The main condition of Theorem 1 on the index $r$ can be similarly understood
from the following lemma:
\begin{Lm}
If $\om\in L^p(\tor)$ for $p>4/3$ and $\bu=\bK*\om$, then for any $h\in {\cal
H}_r$ it holds that $\bu h(\om)\in L^1(\tor)$
when $r={{3}\over{2}}p-1$ for ${{4}\over{3}}<p<2$, $r<p$ for $p=2$, and $r=p$
for $p>2$.
\end{Lm}
{\it Remark:} For convenience in the proof below, and in all later proofs,
we employ an equivalent definition of the class of functions
\be
{\cal H}_p := \left\{ h |\,\,\,h\in C^1({\Bbb R}),
\,\,|h'(\omega)|\leq C |\omega|^{p-1} \,\,{\rm for} \,\,|\om|\geq R
\,\,{\rm for} \,\,{\rm some} \,\, C,R>0 \right\} \lb{Cp1_new} \ee
We will make the argument then assuming that $R=0$ so that the bound
in (\ref{Cp1_new}) above holds globally. In fact, when $R>0$ it is easy
to bound the contributions from the small-$\om$ regions of integration
over space and time by terms proportional to $\|h'\|_{L^\infty[-R,R]},
\|h''\|_{L^\infty[-R,R]}$, assuming that the latter are finite. So we
lose no generality and simplify the arguments by taking $R=0$.

{\it Proof of Lemma:}
We first note the definition $\bK:=\grad^\perp G$ where $G$ is the Greens
function
of $-\bigtriangleup$ on $\tor$.
Then $\bu\in W^{1,p}(\tor)$ because $\|\bu\|_p\leq \|\bK\|_1 \|\om\|_p$ by
Young's inequality and
$\|\grad\bu\|_p\leq C\|\om\|_p$ by the Calder\'{o}n-Zygmund inequality. Hence,
by the same Sobolev imbedding as before,
$\bu\in L^{p'}(\tor)$ for ${{1}\over{p'}}={{1}\over{p}}-{{1}\over{2}}$ when
${{4}\over{3}}<p<2$ and for any {\it finite} $p'\geq 1$ when $p=2$,
and for $p'=\infty$ when $p>2$. Then, by definition of ${\cal H}_r$,
\be \|\bu h(\om)\|_1 \leq ({\rm const.})\|\bu|\om|^{r}\|_1\leq ({\rm
const.})\|\bu\|_{p'}\|\om\|^{r}_{rq'} \lb{uh-bd-0} \ee
with ${{1}\over{q'}}=1-{{1}\over{p'}}$. When ${{4}\over{3}}<p<2$, then
${{1}\over{q'}}={{3}\over{2}}-{{1}\over{p}}$
and $rq'=p$ for $r={{3}\over{2}}p-1.$ On the other hand, when $p>2$,
then $q'=1$, and $rq'=p$ for $r=p$. Lastly, in the critical case
$p=2$, the only requirement is that $q'>1$. Then $rq'\leq p$ can be
satisfied for any $r<p$ by an appropriate choice of $q'>1$.
Thus, for the given definitions of $r$,
\be \|\bu h(\om)\|_1 \leq ({\rm const.})\|\om\|^{r+1}_{p} \lb{uh-bd} \ee
because $\|\bu\|_{p'}\leq C\|\bu\|_{W^{1,p}}\leq C'\|\om\|_p.$ $\Box$

{\it Proof of Theorem 1:} We consider the filtered balance equation
(\ref{filt-bal}):
$$  \partial_th(\om_\en) + \grad\bdot[\bu_\en h(\om_\en) +
h'(\om_\en)\bsigma_\en] =
h''(\om_\en)\grad\om_\en\bdot \bsigma_\en. $$
and, just as in \cite{DR00}, we show that every term on the lefthand side has a
limit in the sense of distributions for $\en\rightarrow 0$.
We show first that $h(\om_\en)\rightarrow h(\om)$. In fact, by the mean-value
theorem, $h(\om_\en)-h(\om)=h'(\bar{\om}_\en)(\om_\en-\om)$
for
$\bar{\om}_\en(\bx,t)=\lambda(\bx,t)\om(\bx,t)+(1-\lambda(\bx,t))
\om_\en(\bx,t)$ with some $0\leq\lambda(\bx,t)\leq 1$. Then, in
the notations of Lemma 1, we have
$$ |h(\om_\en(\bx,t))-h(\om(\bx,t))|\leq ({\rm const.})
|\bar{\om}_\en(\bx,t)|^{r-1}|\om_\en(\bx,t)-\om(\bx,t)| $$
and thus by H\"{o}lder inequality
\begin{eqnarray}
\|h(\om_\en(t))-h(\om(t))\|_{q'}  & \leq &  ({\rm const.})
\|\bar{\om}_\en(t)\|^{r-1}_{rq'} \|\om_\en(t)-\om(t)\|_{rq'} \cr
                             \,   & \leq &  ({\rm const.}) \|\om(t)\|^{r-1}_p
\|\om_\en(t)-\om(t)\|_p. \lb{h-Lq'-bd}
\end{eqnarray}
By the properties of the mollifier, $\lim_{\en\rightarrow 0}
\|\om_\en(t)-\om(t)\|_p=0$ for a.e. $t\in [0,T]$, and thus
$\lim_{\en\rightarrow 0}\|h(\om_\en(t))-h(\om(t))\|_{q'}=0$. To complete the
argument, we use the uniform bound
\be \|h(\om_\en(t))-h(\om(t))\|_{q'} \leq ({\rm const.})
\|\om\|^{r}_{L^\infty(0,T;L^p(\tor))} \lb{h-Linf-bd} \ee
to conclude by dominated convergence that $\lim_{\en\rightarrow
0}\|h(\om_\en)-h(\om)\|_{L^{q'}([0,T]\times \tor)}=0$, which
implies convergence $h(\om_\en)\rightarrow h(\om)$ in sense of distributions.

We show next for the middle term that $\bu_\en h(\om_\en)\rightarrow \bu
h(\om)$. In fact, with notations again as in Lemma 1,
\begin{eqnarray}
\, &  & \|\bu_\en(t) h(\om_\en(t))-\bu(t) h(\om(t))\|_1  \leq
\|\bu_\en(t)-\bu(t)\|_{p'}\|h(\om_\en(t))\|_{q'}+\|\bu(t)\|_{p'}
\|h(\om_\en(t))-h(\om(t))\|_{q'} \cr
\, &  & \,\,\,\,\,\,\,\,\,\,\,\,\,\,\,\,\,\,\,\,\,\,\,\,\,\,
        \leq ({\rm const.})\|\bu_\en(t)-\bu(t)\|_{p'}\|\om(t)\|_p^{r}
+\|\bu(t)\|_{p'}\|h(\om_\en(t))-h(\om(t))\|_{q'}.
\lb{uh-L1-bd}
\end{eqnarray}
Thus, we see that $\lim_{\en\rightarrow 0}\|\bu_\en h(\om_\en(t))-\bu
h(\om(t))\|_1=0$ for a.e. $t\in [0,T]$. In this case we have the
uniform bound
\be \|\bu_\en h(\om_\en(t))-\bu h(\om(t))\|_1 \leq ({\rm
const.})\|\om\|^{r+1}_{L^\infty(0,T;L^p(\tor))} \lb{uh-Linf-bd} \ee
so that we can use Lebesgue's theorem again to infer $\lim_{\en\rightarrow
0}\|\bu_\en h(\om_\en)-\bu h(\om)\|_{L^1([0,T]\times \tor)}=0$,
which gives the result.

Finally, for the third term we show that $h'(\om_\en)\bsigma_\en\rightarrow
\bzed$ as a distribution. We use the definition
$\bsigma_\en=(\bu\om)_\en-\bu_\en\om_\en$ and the H\"{o}lder inequality
\be \|h'(\om_\en(t))\bsigma_\en(t)\|_1  \leq
\|h'(\om_\en(t))\|_{p/(r-1)}
\|(\bu(t)\om(t))_\en-\bu_\en(t)\om_\en(t)\|_{p/(p-r+1)} \lb{hs-Hold} \ee
along with $\|h'(\om_\en(t))\|_{p/(r-1)}\leq ({\rm
const.})\|\om_\en(t)\|_p^{r-1}$ and the triangle inequality
\begin{eqnarray}
\,  &  &  \|(\bu(t)\om(t))_\en-\bu_\en(t)\om_\en(t)\|_{p/(p-r+1)} \leq
\|(\bu(t)\om(t))_\en-\bu(t)\om(t)\|_{p/(p-r+1)}   \cr
\,  &  &
\,\,\,\,\,\,\,\,\,\,\,\,\,\,\,\,\,\,\,\,\,\,\,\,\,\,
\,\,\,\,\,\,\,\,\,\,\,\,\,\,\,\,\,\,\,\,\,\,\,\,\,\,
         +\|\bu(t)\|_{p'}\|\om(t)-\om_\en(t)\|_{p} +
\|\bu(t)-\bu_\en(t)\|_{p'}\|\om_\en(t)\|_p \lb{sigma-bd}
\end{eqnarray}
to infer that $\lim_{\en\rightarrow
0}\|(\bu(t)\om(t))_\en-\bu_\en(t)\om_\en(t)\|_{p/(p-r+1)}=0$ for a.e. $t\in
[0,T]$. Note that we
have used $\|\bu(t)\om(t)\|_{p/(p-r+1)}\leq \|\bu(t)\|_{p/(p-r)}\|\om(t)\|_p$
and $(p-r)/p<p'$. Again a uniform bound on
$\|h'(\om_\en(t))\bsigma_\en(t)\|_1$ like that in (\ref{uh-Linf-bd}) completes
the argument. Gathering these results, we see
that the entire lefthand side of (\ref{filt-bal}) approaches $\partial_th(\om)
+ \grad\bdot[\bu h(\om)]$ in the sense of distributions
as $\en\rightarrow 0$. Obviously this limit is independent of the mollifier
$\varphi$ and the righthand side $-Z_{h,\en}(\om)
= h''(\om_\en)\grad\om_\en\bdot \bsigma_\en$ has the same limit. This gives the
first half of Theorem 1.

The second half of the theorem for the particular choice
$h(\om)={{1}\over{2}}|\om|^2$ follows by the same argument as in \cite{DR00}.
In this proof, the balance (\ref{filt-bal}) is replaced by
\be
\partial_t({{1}\over{2}}\om\om_\en)+\grad\bdot[({{1}\over{2}}\om\om_\en)\bu] =
-\widetilde{Z}_\en(\om) \lb{DR-bal} \ee
where an easy calculation gives
\be \widetilde{Z}_\en(\om) :=
{{1}\over{2}}\om\grad\bdot[(\om\bu)_\en]
-{{1}\over{2}}\om(\bu\bdot\grad)\om_\en. \lb{Ztild} \ee
An argument exactly like the previous one shows that, when $p>2,$ the
distributional limit $\lim_{\en\rightarrow 0}\widetilde{Z}_\en(\om)$
exists and equals
$-\partial_t({{1}\over{2}}|\om|^2)-\grad\bdot[({{1}\over{2}}|\om|^2)\bu]=
Z(\om)$. In addition, a simple calculation
using the incompressibility of the velocity field shows that the expression
appearing in (\ref{4/5law}) in Theorem 1 can be written
\be \int d^2\bl \,\,\grad\varphi_\en(\bl)\bdot\Delta_\bl\bu|\Delta_\bl\om|^2 =
\grad\bdot[\bu(\om^2)_\en-(\bu\om^2)_\en]
+4\widetilde{Z}_\en(\om). \lb{alt-exp} \ee
As before, it is easy to show for $p>2$ that
$\bu(\om^2)_\en-(\bu\om^2)_\en\rightarrow \bzed$ as a distribution when
$\en\rightarrow 0$.
Hence, it follows that the limits of ${{1}\over{4}}\int d^2\bl
\,\,\grad\varphi_\en(\bl)\bdot\Delta_\bl\bu|\Delta_\bl\om|^2$ and
$\widetilde{Z}_\en(\om)$ are also the same. That proves the second half of
Theorem 1. $\Box$


\noindent {\bf 2.2. Proof of Theorem 2}

\noindent A result on global conservation corresponding to the local result in
Theorem 2 was already proved in \cite{Ey96} but with an additional smoothness
assumption that $\om\in L^p(0,T;B^{s,\infty}_p(\tor))$. Here we show that
conservation holds without any such a smoothness assumption. Let
$\tau_\en(f,g):
=(fg)_\en-f_\en g_\en$ where $f_\en= \varphi_\en*f$. Then, we make
use of the following key estimate:
\begin{Lm}
Let $\om\in L^p(\tor)$ and $\bu\in W^{1,p}(\tor)$ for $p\geq 2$, and let
$\grad\bdot\bu=0$. Then
\be \|\grad\bdot\tau_\en(\bu,\om)\|_{L^{p/2}}\leq C\|\bu\|_{W^{1,p}}
                                                     \|\om\|_{L^p} \lb{key} \ee
with a constant $C$ independent of $\en$.
\end{Lm}
\noindent {\it Proof:} Note that
\be \grad\bdot\tau_\en(\bu,\om)=\grad\bdot[(\bu\om)_\en-\bu\om_\en]+
                              (\bu-\bu_\en)\bdot\grad\om_\en.   \lb{gradtau}
\ee
The first term is handled in exactly the same manner as in Lemma II.1 of
\cite{DP-L}.
However, it is easy to see that
\be \|\bu-\bu_\en\|_{L^p}\leq \en\|\grad\bu\|_{L^p} \leq
                              \en\|\bu\|_{W^{1,p}} \lb{ineq1} \ee
and
\be \|\grad\om_\en\|_{L^p} \leq \en^{-1}\|\grad\varphi\|_{L^1}
                                        \|\om\|_{L^p}.   \lb{ineq2} \ee
These control the second term. $\Box$
\begin{Cor}
Under the same hypotheses, let $r_\en:= -\grad\bdot\tau_\en(\bu,\om)$. Then
$\lim_{\en\rightarrow 0}r_\en=0$ strong in $L^{p/2}(\tor)$ for $p\geq 2$.
\end{Cor}
\noindent {\it Proof:} Since $\lim_{\en\rightarrow 0}r_\en=0$ for smooth
$\bu,\om$,
one can obtain the result for all $\om\in L^p(\tor),\bu\in W^{1,p}(\tor)$ by
an approximation argument using the estimate in Proposition 1. $\Box$

If $\bu$ is related to $\om$ by the Biot-Savart formula, $\bu=\bK*\om$, then
$\tau_\en(\bu,\om)=\bsigma_\en$ in the earlier notation. In particular, we see
that
\be \partial_t\om_\en +(\bu_\en\bdot\grad)\om_\en = r_\en \lb{filtereq} \ee
for a weak Euler solution.

\noindent {\it Proof of Theorem 2:} Using (\ref{filtereq}) we get
\be \partial_th(\om_\en) +(\bu_\en\bdot\grad)h(\om_\en) =
                                       h'(\om_\en)r_\en. \lb{filtbal} \ee
It was proved in Theorem 1 that
\be \partial_th(\om_\en) +(\bu_\en\bdot\grad)h(\om_\en) \longrightarrow
     \partial_th(\om) +(\bu\bdot\grad)h(\om) \lb{lhs} \ee
in the sense of distributions for all such $h$. Furthermore, for any
$h\in C^1$ with $h'\in L^\infty$,
\be \|h'(\om_\en)r_\en\|_{L^1}\rightarrow 0. \lb{rhs} \ee
Having proved that (\ref{loc-cons}) holds for $h$ with $h'\in L^\infty$
we then extend it to the general $h$ in the theorem statement by an
approximation argument, as in Corollaries II.1-2 in \cite{DP-L}. $\Box$

\noindent {\it Remark:} The smoothness assumed in the earlier proof of
\cite{Ey96} is not necessary to obtain conservation, but only to provide
an estimate of the rate of the vanishing of the flux. With the assumption
that $\om\in L^p(0,T;B^{s,\infty}_p(\tor))$ the bounds above can be
improved as follows. General estimates in Besov spaces give
$\|\grad\om_\en(t)\|_p\leq C\en^{s-1}\|\om(t)\|_{B^{s,\infty}_p}$
and $\sup_{|\bl|<\en}\|\Delta_\bl\om(t)\|_p\leq \en^s
\|\om(t)\|_{B^{s,\infty}_p}$. See \cite{Treib}, or Appendix C of \cite{Ey96}.
Just as in \cite{Ey96} this gives
\be  \|Z_{h,\en}(\om(t))\|_1 \leq ({\rm const.}) \en^{2s} \|\om(t)\|_p^{r-1}
\|\om(t)\|_{B^{s,\infty}_p}^2
                            \leq ({\rm const.}) \en^{2s}
\|\om(t)\|_{B^{s,\infty}_p}^{r+1} \lb{Zt-fin-bd-B} \ee
Because $r+1<p$, integrating over $t\in [0,T]$ gives
\be \|Z_{h,\en}(\om)\|_{L^1([0,T]\times\tor)} \leq ({\rm const.}) \en^{2s}
\|\om\|_{L^p(0,T;B^{s,\infty}_p(\tor))}^{r+1}. \lb{fin-Z-bd-B} \ee
Thus, $\lim_{\en\rightarrow 0}Z_{h,\en}(\om)=0$ as before, but with an
estimate of the rate. The bound $O(\en^{2s})$ is in agreement with the
estimate given by the heuristic argument in the Introduction.

\vspace{.2in}

\noindent {\bf 2.3. Proof of Theorem 3}

\noindent Theorem 3 is a consequence of the following technical lemma:
\begin{Prop}
Consider a sequence $\{\omega^\en|\en>0\}$ and $\bu^\en=\bK*\om^\en$ given
by the Biot-Savart formula, with the following properties:
\be \sup_{\en>0}\|\om^\en\|_{L^r(0,T;B^{0,\infty}_p(\tor))}<\infty
    \lb{besov} \ee
for $r,p\in [2,\infty]$, and
\be \sup_{\en>0}\|\bu^\en\|_{{\rm Lip}(0,T;H^{-L}(\tor))} <\infty
     \lb{lipschitz} \ee
for some $L>3$. Then, there exist $\om$ and $\bu=\bK*\om$ with
\be \om \in L^r(0,T;B^{0,\infty}_p(\tor))
    \lb{omega-space} \ee
and
\be \bu\in {\rm Lip}(0,T;H^{-L}(\tor)). \lb{velocity-space} \ee
and, furthermore, there exists a subsequence of $\om^\en,\bu^\en$ along which
\be \om^\en\rightharpoonup \om \,\,\,\, weak-* \,\,\,\, {\rm in} \,\,\,\,
      L^r(0,T;B^{0,\infty}_p(\tor)), \lb{om-weak-*} \ee
\be \om^\en\rightarrow \om \,\,\,\, strong \,\,\,\, in \,\,\,\,
      L^r(0,T;W^{-s,q}(\tor)), \lb{om-strong} \ee
for some $q>p$ and $s>2\left({{1}\over{p}}-{{1}\over{q}}\right)$ and
for $t=\min\{r,q\}\geq 2$,
\be \bu^\en\rightarrow \bu  \,\,\,\, strong \,\,\,\, in \,\,\,\,
      L^t([0,T]\times \tor). \lb{u-strong} \ee
\end{Prop}

\noindent {\it Proof:} The first statement (\ref{om-weak-*}) on weak-*
convergence
of $\om^\en$ to $\om\in L^r(0,T;B^{0,\infty}_p(\tor))$ is a simple consequence
of the Banach-Alaoglu theorem.

We derive the second statement from the Aubin-Lions compactness criterion
(see \cite{Lions}, Theorem 5.1 or \cite{Temam}, Theorem III.2.1). Note first
that there is the continuous embedding $B^{0,\infty}_p(\tor)\subset
B^{-s',q}_q(\tor)
=W^{-s',q}(\tor)$ for each $q>p$ and
$s'>2\left({{1}\over{p}}-{{1}\over{q}}\right)$
(see \cite{Treib}, Theorem 2.7.1 and Prop.2.3.2/2). Therefore, from
(\ref{besov}),
\be \sup_{\en>0}\|\om^\en\|_{L^r(0,T;W^{-s',q}(\tor))}<\infty
    \lb{sobolev} \ee
On the other hand, from (\ref{lipschitz}),
\be \sup_{\en>0}\left\|{{d\om^\en}\over{dt}}\right\|_{L^\infty(0,T;
H^{-(L+1)}(\tor))}. \lb{om-lipschitz} \ee
Furthermore, for $s>s'$ and $L+1>s$ there are continuous embeddings
\be W^{-s',q}(\tor))\subset W^{-s,q}(\tor))\subset H^{-(L+1)}(\tor),
\lb{embeds} \ee
and the first embedding is compact by the Rellich-Kondrachov theorem (see \cite
{Mazja}, Chapter 12). Hence, we conclude that $\{\omega^\en|\en>0\}$ is compact
in $L^r(0,T;W^{-s,q}(\tor))$ and contains a strongly convergent subsequence.

To obtain the third result we remark that one may choose $0<s<1$ and that the
mapping
$\om \mapsto\bu = \bK*\om$ is continuous from $W^{-s,q}(\tor)$ into
$W^{1-s,q}(\tor)$,
because of the continuity of the singular integral operator ${\bf T}(\om)=
(\grad\bK)*\om$ from $W^{-s,q}(\tor)$ into itself (for example, see
\cite{Torres},
Theorem 3.2.1) and the bound $\|\bu\|_{W^{1-s,q}(\tor)}\leq ({\rm const.})
\left[\|\bu\|_{W^{-s,q}(\tor)}+\|\grad\bu\|_{W^{-s,q}(\tor)}\right]$ (see
\cite{Treib},
Theorem 2.3.8). Of course, convergence of $\bu^\en\rightarrow\bu$ strong
in $L^r(0,T;W^{1-s,q}(\tor))$ implies at once convergence strong in
$L^t([0,T]\times \tor).$ $\Box$

\noindent {\it Proof of Theorem 3:} The proof is very straightforward and quite
similar to that of DiPerna and Majda in \cite{DP-M} for $\om_0\in L^p$ with
$p\geq 2$
(the easier case than $1<p<2$). In fact, the Lipschitz estimate in time
\be \sup_{\nu>0}\|\bu^\nu\|_{{\rm Lip}(0,T;H^{-L}(\tor))}<\infty
\lb{nu-lipschitz} \ee
holds for the $2D$ Leray solutions with initial energy finite,
$E_0:={{1}\over{2}}
\|\bu_0\|_{L^2}^2<\infty$, which is part of our assumption. See Section 2A and
Appendix A of \cite{DP-M}, for example. But, in that case, the Proposition 1
applies,
with $r=p=2$. The limiting velocity $\bu$ is easily seen to be a weak solution
of the 2D Euler equation in the velocity-pressure formulation, because of the
third
result (\ref{u-strong}), the strong $L^t$ convergence $\bu^\nu\rightarrow \bu$
with $t>2$. Obviously, more general versions of Theorem 3 for any $r,p\in
[2,\infty]$
could be proved, with (\ref{besov}) replacing (\ref{om-besov}) in the
hypothesis.

The statements on the weak solutions in the vorticity-velocity formulation
follow from arguments very similar to those earlier in Lemma 1, but now using
the density of $C^\infty(\tor)$ in $B^{0,1}_q(\tor)$ (\cite{Treib}, Theorem
2.3.3).
So, we just verify the required regularity of $\bu$. Because of the hypothesis
on $\om$, $\bu\in L^r(0,T;B^{1,\infty}_p(\tor))$. Indeed, using the
Calder\'{o}n-Zygmund inequality it is easy to show that $\grad\bu\in
L^r(0,T;B^{0,\infty}_p(\tor))$, and this is equivalent to the first statement
(\cite{Treib}, Theorem 2.3.8). Then, for any $p''<p'$ with $p'$ defined by
${{1}\over{p'}}={{1}\over{p}}-{{1}\over{2}}$, one has the continuous embedding
$B^{1,\infty}_p(\tor)\subset B^{0,1}_{p''}(\tor)$. In fact,
$B^{1,\infty}_p(\tor)
\subset B^{\en,\infty}_{p''}(\tor)$ for ${{1}\over{p''}}:={{1}\over{p'}}+
{{\en}\over{2}}$ for any small $\en>0$ (\cite{Treib}, Theorem 2.7.1),
but then $B^{\en,\infty}_{p''}(\tor)\subset B^{0,1}_{p''}(\tor)$
by an elementary imbedding (\cite{Treib}, Prop. 2.3.2/2). Now, precisely
for $p>{{4}\over{3}}$, one has ${{1}\over{p'}}<{{1}\over{q}}$ with
${{1}\over{p}}+{{1}\over{q}}=1$. Thus, it is possible to choose $p''>q$
but still $p''<p'$. Also, with $t$ defined by ${{1}\over{r}}+{{1}\over{t}}=1$,
$r\geq t$ for $r\geq 2$. In that case, $\bu\in L^t(0,T;B^{0,1}_q(\tor))$,
whereas $L^r(0,T;B^{0,\infty}_p(\tor))=[L^t(0,T;B^{0,1}_q(\tor))]^*$,
the Banach dual (\cite{Treib}, Theorem 2.11.2). To conclude the proof,
we just note that, if $\psi\in C_0^\infty([0,T]\times\tor)$, then also
$(\partial_t+\bu\bdot\grad)\psi\in L^t(0,T;B^{0,1}_q(\tor))$
(see \cite{Treib}, Lemma 3.3.1). $\Box$

\newpage

\noindent {\bf Acknowledgements} I would like to thank the Institute of
Theoretical Physics at Santa Barbara
for their kind hospitality in the spring of 2000. This work was begun there
during my stay for the program
on ``Physics  of Hydrodynamic Turbulence.'' It was completed during my
sabbatical at Johns Hopkins
University in autumn 2000, and I would like to acknowledge my hosts there,
S. Chen and C. Meneveau. I also wish to thank W. E, U. Frisch,
C. D. Levermore, P.-L. Lions, and E. Vanden-Eijnden for very useful
conversations and correspondence.


\begin{thebibliography}{99}
\bibitem[1]{Kr67} R. H. Kraichnan, ``Inertial ranges in two-dimensional
turbulence,'' Phys. Fluids {\bf 10} 1417-1423 (1967)
\bibitem[2]{Leith} C. E. Leith, ``Diffusion approximation for two-dimensional
turbulence,'' Phys. Fluids
{\bf 11} 671-672 (1968)
\bibitem[3]{Ba69} G. K. Batchelor, ``Computation of the energy spectrum in
homogeneous two-dimensional
turbulence,'' Phys. Fluids Suppl. II {\bf 12} 233-239 (1969)
\bibitem[4]{On49} L. Onsager, ``Statistical hydrodynamics,'' Nuovo Cim. Suppl.
{\bf 6} 279-289 (1949)
\bibitem[5]{DR00} J. Duchon and R. Robert, ``Inertial energy dissipation for
weak solutions of incompressible
Euler and Navier-Stokes equations,'' Nonlinearity {\bf 13} 249-255 (2000)
\bibitem[6]{Poly} A. M. Polyakov, ``The theory of turbulence in two
dimensions,'' Nucl. Phys. B {\bf 396} 367-385 (1993)
\bibitem[7]{Ey96} G. L. Eyink, ``Exact results for stationary turbulence in 2D:
consequences of vorticity
conservation,'' Physica D {\bf 91} 97-142 (1996)
\bibitem[8]{C-W} P. Constantin and J. Wu, ``The inviscid limit for nonsmooth
vorticity,'' Indiana Univ. Math. J.
{\bf 45} 67-81 (1996)
\bibitem[9]{MK00} C. Meneveau and J. Katz, ``Scale-invariance and turbulence
models for large-eddy simulation,''
Annu. Rev. Fluid Mech. {\bf 32} 1-32 (2000)
\bibitem[10]{BMPS} M. E. Brachet, M. Meneguzzi, H. Politano, and P. L. Sulem,
``The dynamics of freely decaying
two-dimensional turbulence,'' J. Fluid Mech. {\bf 194} 333-49 (1988)
\bibitem[11]{DP-M} R. J. DiPerna and A. J. Majda, ``Concentrations in
regularizations for 2-D
incompressible flow,'' Commun. Pure Appl. Math. {\bf XL} 301-345 (1987)
\bibitem[12]{Treib} H. Treibel, {\it Theory of Function Spaces}
(Birkh\"{a}user, Basel, 1983)
\bibitem[13]{CET} P. Constantin, W. E, and E. S. Titi, ``Onsager's conjecture
on the energy conservation
for solutions of Euler's equation,'' Commun. Math. Phys. {\bf 165} 207-209
(1994)
\bibitem[14]{Lions96}P.-L. Lions, {\it Mathematical Topics in Fluid Mechanics.
Vol. 1. Incompressible Models}. (Clarendon Press, Oxford, 1996)
\bibitem[15]{DP-L} R. J. DiPerna and P. L. Lions, ``Ordinary differential
equations, transport theory, and Sobolev spaces,'' Invent. Math. {\bf 98}
511-547 (1989).
\bibitem[16]{E-vdE2} W. E and E. Vanden-Eijnden, ``Remarks on 2D turbulence,''
preprint.
\bibitem[17]{GKB} D. Bernard, K. Gaw\c{e}dzki, and A. Kupiainen, ``Slow modes
in passive advection,''
J. Stat. Phys. {\bf 90} 519-569 (1998)
\bibitem[18]{LJ-R} Y. Le Jan and O. Raimond, ``Integration of Brownian vector
fields,'' preprint
{\it math.PR/9909147}
\bibitem[19]{E-vdE} W. E and E. Vanden-Eijnden, ``Generalized flows, intrinsic
stochasticity, and
turbulent transport,'' preprint {\it nlin.CD/0003028}.
\bibitem[20]{Shnir} A. Shnirelman, ``Weak solutions with decreasing energy of
incompressible Euler equations,''
Commun. Math. Phys. {\bf 210} 541-603 (2000)
\bibitem[21]{Yud} V. I. Yudovich, ``Non-stationary flow of an ideal
incompressible liquid,'' Zh. Vych.
Mat. {\bf 3} 1032-1066 (1963)
\bibitem[22]{Const-spect}P. Constantin, ``The Littlewood-Paley spectrum in
two-dimensional turbulence," Theor. Comp. Fluid Dyn. {\bf 9} 183-9 (1997)
\bibitem[23]{Kr71} R. H. Kraichnan, ``Inertial-range transfer in two- and
three-dimensional turbulence,'' J. Fluid Mech. {\bf 47} 525-535 (1971)
\bibitem[24]{Ey95} G. L. Eyink, ``Besov spaces and the multifractal
hypothesis,''
J. Stat. Phys. {\bf 78} 353-375 (1995).
\bibitem[25]{BC94}H. Bahouri and J.-Y. Chemin, ``Equations de transport
relatives \`{a} des champs
de vecteurs non-lipschitziens et m\'{e}canique des fluides,'' Arch. Rat. Mech.
Anal. {\bf 127}
159-181 (1994)
\bibitem[26]{Kr75}R. H. Kraichnan, ``Statistical dynamics of two-dimensional
flow,'' J. Fluid Mech.
{\bf 67} 155-175 (1975)
\bibitem[27]{Kr74}R. H. Kraichnan, ``Convection of a passive scalar by a
quasi-uniform random straining field,'' J. Fluid Mech.
{\bf 64} 737-762 (1974)
\bibitem[28]{Lax}P. D. Lax, ``Shock waves and entropy,'' in: {\it Contributions
to Nonlinear Functional Analysis.}
E. H. Zarantonello, Ed. (Academic Press, New York, 1971)
\bibitem[29]{Vishik} M. Vishik, ``Incompressible flows of an ideal fluid
with vorticity in borderline spaces of Besov type,'' Ann. de l'Ecole
Norm. Sup. {\bf 32} 769-812 (1999)
\bibitem[30]{Lions} J. L. Lions, {\it Quelques M\'ethodes de Resolution de
Probl\`emes
aux Limites Non Lineaires}. (Dunot Gauthier-Villars, Paris, 1969)
\bibitem[31]{Temam}R. Temam, {\it Navier-Stokes Equations: Theory and Numerical
Analysis.} (Elsevier, Amsterdam, 1984)
\bibitem[32]{Mazja}V. G. Maz'ja, {\it Sobolev spaces}. (Springer-Verlag,
Berlin, 1985)
\bibitem[33]{Torres}R. H. Torres, {\it Boundedness Results for Operators with
Singular Kernels on Distributions}. (American Mathematical Society,
Providence, RI, 1991)

\end{thebibliography}
\end{document}